\documentclass[a4paper]{amsproc}

\usepackage{amssymb}
\usepackage{pstricks}
\usepackage{pstcol}
\usepackage{graphicx}
\usepackage{amsmath}
\usepackage{amsfonts}
\usepackage{latexsym}

\theoremstyle{plain}

\theoremstyle{definition}
\newtheorem{exm}{Example}[section]

\numberwithin{equation}{section}

\title[Logarithmic integrals]{The integrals in Gradshteyn and Ryzhik.  \\
Part 12: Some logarithmic integrals}

\subjclass[2000]{Primary 33}


\author[Victor H. Moll]{Victor H. Moll}
\address{Department of Mathematics,
Tulane University, New Orleans, LA 70118}

\author[Ronald A. Posey]{Ronald A. Posey}
\address{Department of Mathematics, Baton Rouge Community College, 
Baton Rouge, LA 70806}

\address{\hfill{\it Received  02 07 2008, revised 16 03 2009}
\newline Departamento de Matem\'atica
\newline
Universidad T\'ecnica Federico Santa Mar\'{\i}a
\newline  Casilla 110-V,
\newline Valpara\'{\i}so, Chile}

\email{vhm@math.tulane.edu}
\email{raposey@gmail.com}

\thanks{The first author wishes to thank the 
partial support of NSF-DMS 0713836}

\begin{document}

{\begin{flushleft}\baselineskip9pt\scriptsize {\bf SCIENTIA}\newline
Series A: {\it Mathematical Sciences}, Vol. 18 (2009), 77-84
\newline Universidad T\'ecnica Federico Santa Mar{\'\i}a
\newline Valpara{\'\i}so, Chile
\newline ISSN 0716-8446
\newline {\copyright\space Universidad T\'ecnica Federico Santa
Mar{\'\i}a\space 2009}
\end{flushleft}}

\vspace{10mm} \setcounter{page}{1} \thispagestyle{empty}

\begin{abstract}
We present the evaluation of some logarithmic integrals. The integrand 
contains a rational function with complex poles. The methods are 
illustrated with 
examples found in the classical table of 
integrals by I. S. Gradshteyn and I. M. Ryzhik.
\end{abstract}

\maketitle

\newcommand{\nn}{\nonumber}
\newcommand{\ba}{\begin{eqnarray}}
\newcommand{\ea}{\end{eqnarray}}
\newcommand{\ift}{\int_{0}^{\infty}}
\newcommand{\ione}{\int_{0}^{1}}
\newcommand{\ifft}{\int_{- \infty}^{\infty}}
\newcommand{\no}{\noindent}
\newcommand{\Ftwo}{{{_{2}F_{1}}}}
\newcommand{\realpart}{\mathop{\rm Re}\nolimits}
\newcommand{\imagpart}{\mathop{\rm Im}\nolimits}

\newtheorem{Definition}{\bf Definition}[section]
\newtheorem{Thm}[Definition]{\bf Theorem} 
\newtheorem{Example}[Definition]{\bf Example} 
\newtheorem{Lem}[Definition]{\bf Lemma} 
\newtheorem{Note}[Definition]{\bf Note} 
\newtheorem{Cor}[Definition]{\bf Corollary} 
\newtheorem{Prop}[Definition]{\bf Proposition} 
\newtheorem{Problem}[Definition]{\bf Problem} 
\numberwithin{equation}{section}

\section{Introduction} \label{intro} 
\setcounter{equation}{0}

The classical table of integrals by I. Gradshteyn and I. M. Ryzhik 
\cite{gr} contains many entries from  the family
\begin{equation}
\int_{0}^{1} R(x) \, \log x \, dx 
\end{equation}
\noindent
where $R$ is a rational function. For instance, the elementary integral 
$\mathbf{4.231.1}$
\begin{equation}
\int_{0}^{1} \frac{\log  x \, dx}{1+x} = - \frac{\pi^{2}}{12},
\end{equation}
\noindent
is evaluated simply by expanding the integrand in a power series.  In 
\cite{moll-gr9}, the first author and collaborators have presented a systematic
study of integrals of the form
\begin{equation}
h_{n,1}(b) = \int_{0}^{b} \frac{\log t \, dt}{(1+t)^{n+1}}, 
\end{equation}
\noindent
as well as the case in which the integrand has a single purely imaginary pole
\begin{equation}
h_{n,2}(a,b) = \int_{0}^{b} \frac{\log t \, dt}{(t^{2}+a^{2})^{n+1}}. 
\end{equation}

The work presented here deals with integrals where the rational part of the 
integrand is allowed to have arbitrary complex poles. 

\section{Evaluations in terms of polylogarithms} \label{sec-first} 
\setcounter{equation}{0}

In this section we describe the evaluation of 
\begin{equation}
g(a) = \int_{0}^{1} \frac{\log x \, dx}{x^{2} - 2ax + 1},
\end{equation}
\noindent
under the assumption that the denominator has non-real roots, that is, 
$a^{2} < 1$.

The first approach to the evaluation of $g(a)$ 
is based on the factorization of the quartic as 
\begin{equation}
x^{2}-2ax+1 = (x+r_{1})(x+r_{2}),
\end{equation}
\noindent
where $r_{1} = -a + i \sqrt{1-a^{2}}$ and $r_{2} = -a - i \sqrt{1-a^{2}}$. The
partial fraction expansion
\begin{equation}
\frac{1}{(x+r_{1})(x+r_{2})} = \frac{1}{r_{2}-r_{1}} \left( \frac{1}{x+r_{1}} 
- \frac{1}{x+r_{2}} \right), 
\end{equation}
\noindent
yields
\begin{equation}
g(a) = \frac{1}{r_{2}-r_{1}} \int_{0}^{1} \frac{\log x \, dx }{x+r_{1}} 
 - \frac{1}{r_{2}-r_{1}} \int_{0}^{1} \frac{\log x \, dx }{x+r_{2}}.
\end{equation}

These  integrals are computed in terms of the {\em dilogarithm} function
defined by
\begin{equation}
\text{PolyLog}[2,x] := - \int_{0}^{x} \frac{\log(1-t)}{t} \, dt.
\end{equation}
\noindent
A direct calculaton shows that
\begin{equation}
\int \frac{\log x \, dx}{x+a} = \log x \log(1 + x/a) + 
\text{PolyLog}[2,-x/a], 
\end{equation}
\noindent
and thus
\begin{equation}
\int_{0}^{1} \frac{\log x \, dx}{x+a} =
\text{PolyLog} \left[2,-\frac{1}{a} \right].
\end{equation}
\noindent
It follows that
\begin{equation}
g(a) = \frac{1}{r_{2}-r_{1}} 
\left( \text{PolyLog} \left[2, - \frac{1}{r_{1}} \right] -
       \text{PolyLog} \left[2, - \frac{1}{r_{2}} \right] \right). 
\end{equation}

Observe that the real integral $g(a)$ appears here expressed in terms of 
the polylogarithm of complex arguments.

\begin{exm}
The case $a = 1/2$ yields
\begin{equation}
\int_{0}^{1} \frac{\log x \, dx}{x^{2}-x + 1} = 
\frac{i}{\sqrt{3}} \left(\text{PolyLog} \left[2,(1+i \sqrt{3})/2 \right] -
       \text{PolyLog} \left[2, (1-i \sqrt{3})/2 \right] \right). 
\label{int-ex1}
\end{equation}
\noindent
The polylogarithm function is evaluated using the representation
\begin{equation}
(1 + i \sqrt{3})/2 = e^{i \pi/3},
\end{equation}
\noindent
to produce
\begin{eqnarray}
\text{PolyLog} \left[2,(1+i \sqrt{3})/2 \right]  & = & 
\sum_{k=1}^{\infty} \frac{\left[ \tfrac{1}{2}(1 + i \sqrt{3}) \right]^{k}}
{k^{2}} 
  =  \sum_{k=1}^{\infty} \frac{e^{i \pi k/3}}{k^{2}}  \nonumber \\
& = & \sum_{k=1}^{\infty} 
\frac{\cos \left( \frac{\pi k }{3} \right) + i \sin \left( \frac{\pi k}{3} 
\right) }{k^{2}}. 
\nonumber
\end{eqnarray}
\noindent
Similarly
\begin{eqnarray}
\text{PolyLog} \left[2,(1-i \sqrt{3})/2 \right]  & = & 
\sum_{k=1}^{\infty} 
\frac{\cos \left( \frac{\pi k }{3} \right) - i \sin \left( \frac{\pi k}{3} 
\right) }{k^{2}},
\nonumber
\end{eqnarray}
and it follows that
\begin{eqnarray}
\int_{0}^{1} \frac{\log x \, dx}{x^{2}-x + 1}  & = & 
\frac{i}{\sqrt{3}} \left(\text{PolyLog} \left[2,(1+i \sqrt{3})/2 \right] -
       \text{PolyLog} \left[2, (1-i \sqrt{3})/2 \right] \right)
\nonumber \\
& = & - \frac{2}{\sqrt{3}} \sum_{k=1}^{\infty} 
\frac{\sin \left( \frac{\pi k}{3} \right) }{k^{2}}. 
\nonumber
\end{eqnarray}
\noindent

The function $\sin(\pi k/3)$ is periodic, with period $6$, and
repeating values $$\tfrac{\sqrt{3}}{2}, \tfrac{\sqrt{3}}{2}, 0, 
- \tfrac{\sqrt{3}}{2}, - \tfrac{\sqrt{3}}{2}, 0.$$
\noindent
Therefore 

\begin{equation}
\sum_{k=1}^{\infty} \frac{\sin \left( \frac{\pi k}{3} \right) }{k^{2}} = 
\frac{\sqrt{3}}{2} \left( 
\sum_{k=0}^{\infty} \frac{1}{(6k+1)^{2}} 
+ \sum_{k=0}^{\infty} \frac{1}{(6k+2)^{2}} 
- \sum_{k=0}^{\infty} \frac{1}{(6k+4)^{2}} 
- \sum_{k=0}^{\infty} \frac{1}{(6k+5)^{2}} \right). \nonumber
\end{equation}

\medskip

To evaluate these sums, recall the series representatin of the {\em polygamma}
function $\psi(x) = \Gamma'(x)/\Gamma(x)$, given by
\begin{equation}
\psi(x) = - \gamma - \frac{1}{x} + \sum_{k=1}^{\infty} \frac{x}{k(x+k)}. 
\end{equation}
\noindent
Differentiation yields
\begin{equation}
\psi'(x) = - \sum_{k=0}^{\infty} \frac{1}{(x+k)^{2}},
\end{equation}
\noindent
and we obtain
\begin{equation}
\sum_{k=0}^{\infty} \frac{1}{(6k+j)^{2}} = \frac{1}{36} 
\sum_{k=0}^{\infty} \frac{1}{(k+\tfrac{j}{6})^{2}}.
\nonumber
\end{equation}
\noindent
This provides the expression
\begin{equation}
\sum_{k=1}^{\infty} \frac{\sin \left( \frac{\pi k}{3} \right) }{k^{2}} = 
\frac{\sqrt{3}}{72} \left( 
\psi' \left( \tfrac{1}{6} \right) +
\psi' \left( \tfrac{2}{6} \right) -
\psi' \left( \tfrac{4}{6} \right) -
\psi' \left( \tfrac{5}{6} \right) \right).
\end{equation}
\noindent
The integral (\ref{int-ex1}) is
\begin{equation}
\int_{0}^{1} \frac{\log x \, dx}{x^{2}-x + 1}  = 
-\frac{1}{36} \left( 
\psi' \left( \tfrac{1}{6} \right) +
\psi' \left( \tfrac{2}{6} \right) -
\psi' \left( \tfrac{4}{6} \right) -
\psi' \left( \tfrac{5}{6} \right) \right).
\label{int-00}
\end{equation}
\noindent
The identities 
\begin{equation}
\psi(1-x) = \psi(x) + \pi \cot \pi x,
\end{equation}
\noindent
and
\begin{equation}
\psi(2x) = \frac{1}{2} \left( \psi(x) + \psi(x+ \tfrac{1}{2}) \right) + \log 2,
\end{equation}
\noindent
produce
\begin{equation}
\psi' \left( \tfrac{1}{6} \right)   =  5  \psi' \left( \tfrac{1}{3} \right) 
- \tfrac{4 \pi^{2}}{3}, \, 
\psi' \left( \tfrac{2}{3} \right)   =  - \psi' \left( \tfrac{1}{3} \right) 
+ \tfrac{4 \pi^{2}}{3}, \, 
\psi' \left( \tfrac{5}{6} \right)   =  - 5 \psi' \left( \tfrac{1}{3} \right) 
+ \tfrac{16 \pi^{2}}{3}. \nonumber 
\end{equation}
\noindent
Replacing in (\ref{int-00}) yields
\begin{equation}
\int_{0}^{1} \frac{\log x \, dx}{x^{2}-x + 1}  = 
\frac{2 \pi^{2}}{9} - \frac{1}{3} \psi' \left( \frac{1}{3} \right). 
\end{equation}
\noindent
This appears as formula $\mathbf{4.233.2}$ in \cite{gr}. \\

\noindent
{\bf Note}. The  method described in the previous example evaluates 
logarithmic
integrals in terms of the {\em Clausen function}
\begin{equation}
\text{Cl}_{2}(x) := \sum_{k=1}^{\infty} \frac{\sin kx}{k^{2}}. 
\end{equation}

\noindent
{\bf Note}. An identical procedure can be used to evaluate the integrals 
$\mathbf{4.233.1}, \, \mathbf{4.233.3}, \, \mathbf{4.233.4}$ in 
\cite{gr} given by

\begin{equation}
\int_{0}^{1} \frac{\log x \, dx}{x^{2}+x + 1}  = 
\frac{4 \pi^{2}}{27} - \frac{2}{9} \psi' \left( \frac{1}{3} \right),
\end{equation}
\begin{equation}
\int_{0}^{1} \frac{x \, \log x \, dx}{x^{2}+x + 1}  = 
- \frac{7 \pi^{2}}{54} + \frac{1}{9} \psi' \left( \frac{1}{3} \right),
\end{equation}
\noindent
and
\begin{equation}
\int_{0}^{1} \frac{x \, \log x \, dx}{x^{2}-x + 1}  = 
\frac{5 \pi^{2}}{36} - \frac{1}{6} \psi' \left( \frac{1}{3} \right),
\end{equation}
\noindent
respectively. 
\end{exm}

\section{An alternative approach} \label{sec-alter} 
\setcounter{equation}{0}

In this section we present an alternative evaluation for the integral
\begin{equation}
g(a) = \int_{0}^{1} \frac{\log x \, dx}{x^{2} - 2ax + 1},
\label{int-ga}
\end{equation}
\noindent
based on the observation that
\begin{equation}
g(a) = \lim\limits_{s \to 0} \frac{d}{ds} \int_{0}^{1} \frac{x^{s} \, dx}
{x^{2}-2ax+1}.
\end{equation}

The proof discussed here is based on the {\em Chebyshev} polynomials of the 
second kind $U_{n}(a)$, defined by 
\begin{equation}
U_{n}(a) = \frac{\sin[ (n+1)t]}{\sin t}, 
\label{trig-exp}
\end{equation}
\noindent
where $a = \cos t$. The relation with the problem at hand comes from the 
generating function
\begin{equation}
\frac{1}{1-2ax+x^{2}} = \sum_{k=0}^{\infty} U_{k}(a)x^{k}. 
\end{equation}
\noindent
This appears as $\mathbf{8.945.2}$ in \cite{gr}. 

Observe that 
\begin{equation}
\int_{0}^{1} \frac{x^{s} \, dx}{x^{2}-2ax+1}   =  
\sum_{k=0}^{\infty} U_{k}(a) \int_{0}^{1} x^{k+s} dx 
  =  \sum_{k=0}^{\infty} \frac{U_{k}(a)}{k+s+1}. \nonumber 
\end{equation}
\noindent
It follows that
\begin{equation}
\int_{0}^{1} \frac{\log x \, dx}{x^{2}-2ax+1} = - \sum_{k=0}^{\infty} 
\frac{U_{k}(a)}{(k+1)^{2}}. 
\end{equation}
\noindent
Replacing the trigonometric expression (\ref{trig-exp})
for the Chebyshev polynomial, it follows that 
\begin{equation}
\int_{0}^{1} \frac{\log x \, dx}{x^{2}-2ax+1} = - \frac{1}{\sin t} 
\sum_{k=0}^{\infty} 
\, \frac{\sin{kt}}{k^{2}} = - \frac{\text{Cl}_{2}(t)}{\sin t}.
\end{equation}
\noindent
This reproduces the representation
discussed in Section \ref{sec-first}. \\

\noindent
{\bf Note}. The methods presented here give the value of (\ref{int-ga}) 
in terms of the dilogarithm function. The classical values 
\begin{equation}
\text{Cl}_{2} \left( \tfrac{\pi}{2} \right) = 
- \text{Cl}_{2} \left( \tfrac{3\pi}{2} \right) = 
\sum_{k=0}^{\infty} 
\frac{(-1)^{k}}{(2k+1)^{2}} = \text{Catalan}, 
\end{equation}
\noindent
are easy to establish. More sophisticated evaluations appear in \cite{srichoi}.
These are given in terms of the Hurwitz zeta function
\begin{equation}
\zeta(s,q) = \sum_{k=0}^{\infty} \frac{1}{(k+q)^{s}}. 
\end{equation}
\noindent
For instance, the reader will find 
\begin{equation}
\text{Cl}_{2} \left( \tfrac{2 \pi}{3} \right) = 
\sqrt{3} \left( \frac{3^{-s}-1}{2} \zeta(2) + 3^{-s} \zeta(2, \tfrac{1}{3})
\right), 
\end{equation}
\noindent
and 
\begin{equation}
\text{Cl}_{2} \left( \tfrac{\pi}{3} \right) = 
\sqrt{3} \left( 
\frac{3^{-s}-1}{2} \zeta(2) + 6^{-s} \left( \zeta(2, \tfrac{1}{6}) + 
\zeta(2, \tfrac{1}{3}) \right) \right). 
\end{equation}

\noindent
{\bf Note}. Integrals of the form 
\begin{equation}
\int_{0}^{1} R(x) \log \log \frac{1}{x} \, dx 
\end{equation}
\noindent
present new challeges. The reader will find some examples in \cite{luis2}. The 
current version of Mathematica is able to evaluate 
\begin{equation}
\int_{0}^{1} \frac{x \, \log \log 1/x }{x^{4}+x^{2}+1} \, dx = 
\frac{\pi}{12 \sqrt{3}} \left( 6 \log 2 - 3 \log 3 + 8 \log \pi - 12 
\log \Gamma( \tfrac{1}{3}) \right), 
\end{equation}
\noindent
but is unable to evaluate 
\begin{equation}
\int_{0}^{1} \frac{x \, \log  \log 1/x }{x^{4} - \sqrt{2}x^{2}+1} \, dx = 
\frac{\pi}{8 \sqrt{2}} \left( 7 \log \pi - 4 \log \sin \tfrac{\pi}{8} - 
8 \log \Gamma( \tfrac{1}{8}) \right). 
\end{equation}

\section{Higher powers of logarithms} \label{sec-higher} 
\setcounter{equation}{0}

The method of the previous sections can be used to evaluate integrals of the 
form
\begin{equation}
\int_{0}^{1} R(x) \log^{p}x \, dx,
\end{equation}
\noindent
where $R$ is a rational function. The ideas are illustrated with the 
verification of formula $\mathbf{4.261.8}$ in \cite{gr}:
\begin{equation}
\int_{0}^{1} \frac{1-x}{1-x^{6}} \, \log^{2}x \, dx = 
\frac{8 \sqrt{3} \pi^{3} + 351 \zeta(3)}{486}. 
\label{last-example}
\end{equation}

Define 
\begin{eqnarray}
J_{1}  =  \int_{0}^{1} \frac{\log^{2}x \, dx}{1+x}, & & 
J_{2}  =  \int_{0}^{1} \frac{\log^{2}x \, dx}{1-x+x^{2}}, \nonumber \\
J_{3}  =  \int_{0}^{1} \frac{x \, \log^{2}x \, dx}{1-x+x^{2}}, & & 
J_{4}  =  \int_{0}^{1} \frac{\log^{2}x \, dx}{1+x+x^{2}}. \nonumber 
\end{eqnarray}
\noindent
The partial fraction decomposition
\begin{equation}
\frac{1-x}{1-x^{6}} = \frac{1}{3}\frac{1}{1+x} + \frac{1}{6}
\frac{1}{1-x+x^{2}} - \frac{1}{3} \frac{x}{1-x+x^{2}} + 
\frac{1}{2} \frac{1}{1+x+x^{2}},  \nonumber
\end{equation}
\noindent
gives
\begin{equation}
\int_{0}^{1} \frac{1-x}{1-x^{6}} \, \log^{2}x \, dx = 
\frac{1}{3} J_{1} + \frac{1}{6}J_{2} - \frac{1}{3}J_{3} + \frac{1}{2}J_{4}.
\end{equation}

\medskip

\noindent
{\em Evaluation of }$J_{1}$. Consider first
\begin{equation}
\int_{0}^{1} \frac{x^{s}}{1+x} \, dx  =  \sum_{k=1}^{\infty} (-1)^{k-1}
\int_{0}^{1} x^{k+s-1} \, dx
 =  \sum_{k=1}^{\infty} \frac{(-1)^{k-1}}{k+s}. \nonumber
\end{equation}
\noindent
Differentiating twice with respect to $s$ gives
\begin{equation}
J_{1}  =  \int_{0}^{1} \frac{\log^{2}x \, dx}{1+x} = 2 \sum_{k=1}^{\infty} 
\frac{(-1)^{k-1}}{k^{3}} = \frac{3}{2} \zeta(3).
\end{equation}

\medskip

\noindent
{\em Evaluations of }$J_{2}$. Integrating the expansion
\begin{equation}
\frac{x^{s} \, dx}{x^{2}-2ax+1} = \sum_{k=0}^{\infty} \frac{U_{k}(a)}{s+k+1},
\end{equation}
\noindent
and differentiating twice with respect to $s$ yields
\begin{equation}
\int_{0}^{1} \frac{x^{s} \, \log^{2}x \, dx }{x^{2}-2ax+1} = 
2 \sum_{k=0}^{\infty} \frac{U_{k}(a)}{(s+k+1)^{3}}. 
\end{equation}
\noindent
The value $s=0$ yields
\begin{equation}
\int_{0}^{1} \frac{\log^{2}x \, dx }{x^{2}-2ax+1} = 
2 \sum_{k=0}^{\infty} \frac{U_{k}(a)}{(k+1)^{3}}. 
\end{equation}
\noindent
We conclude that
\begin{equation}
J_{2} = 2 \sum_{k=0}^{\infty} \frac{U_{k}(\tfrac{1}{2})}{(k+1)^{3}}. 
\end{equation}
\noindent
The sequence $U_{k}(\tfrac{1}{2})$ is periodic of period $6$ and 
values $1, \, 0, \, -1, \, -1, \, 0, \, 1$. Therefore
\begin{equation}
J_{2} =  
2 \sum_{k=1}^{\infty} \frac{1}{(6k+1)^{3}} -
2 \sum_{k=1}^{\infty} \frac{1}{(6k+3)^{3}} -
2 \sum_{k=1}^{\infty} \frac{1}{(6k+4)^{3}} +
2 \sum_{k=1}^{\infty} \frac{1}{(6k+5)^{3}}. 
\end{equation}
\noindent
This can be written as 
\begin{equation}
J_{2} =  \frac{1}{108} \left( 
\sum_{k=1}^{\infty} \frac{1}{(k+1/6)^{3}} -
\sum_{k=1}^{\infty} \frac{1}{(k+1/2)^{3}} -
\sum_{k=1}^{\infty} \frac{1}{(k+2/3)^{3}} +
\sum_{k=1}^{\infty} \frac{1}{(k+5/6)^{3}} \right). 
\nonumber
\end{equation}
\noindent
Proceeding along the same lines of the previous argument, employing now the 
second derivative of the polygamma function yields
\begin{equation}
J_{2} =  \frac{10 \pi^{3}}{81 \sqrt{3}}. 
\end{equation}

The same type of analysis gives 
\begin{eqnarray}
J_{3} & = & \int_{0}^{1} \frac{x \, \log^{2}x \, dx}{1-x+x^{2}} 
= \frac{5 \pi^{3}}{81 \sqrt{3}} - \frac{2 \zeta(3)}{3}, \nonumber \\
J_{4} & = & \int_{0}^{1} \frac{\log^{2}x \, dx}{1+x+x^{2}} = 
\frac{81 \pi^{3}}{81 \sqrt{3}}. \nonumber 
\end{eqnarray}

\noindent
This completes the proof of (\ref{last-example}). \\

\medskip

The reader is invited to use the method developed here to verify 
\begin{equation}
\int_{0}^{1} \frac{1-x}{1-x^{6}} \, \log^{4}x \, dx = 
\frac{32 \sqrt{3} \pi^{5} + 16335 \zeta(5)}{1458},
\end{equation}
\noindent
and
\begin{equation}
\int_{0}^{1} \frac{1-x}{1-x^{6}} \log^{6}x \, dx = 
\frac{7(256 \sqrt{3} \pi^{7} + 1327995 \zeta(7))}{26244}.
\end{equation}
\noindent
Mathematica 6.2 is capable of producing these results. 

\medskip

The methods discussed here constitute the most elementary approach to 
the evaluations of logarithmic integrals. M. Coffey \cite{coffey-tucson} 
presents some of the more advanced techniques required for the computation
of integrals of the form
\begin{equation}
\int_{0}^{1} R(x) \, \log^{s} x \, dx 
\end{equation}
\noindent
for $s$ real and $R$ a rational function. 

\end{document}